\documentclass[12pt]{amsart}
\usepackage{amssymb,eepic,epic}
\newtheorem{theorem}{Theorem}[section]

\newtheorem{proposition}[theorem]{Proposition}

\theoremstyle{definition}    
\newtheorem{definition}[theorem]{Definition}
\theoremstyle{remark}
\newtheorem{remark}[theorem]{Remark}
\newtheorem{example}[theorem]{Example}
\newcommand\A{\mathcal{A}}
\newcommand\M{\mathcal{M}}
\newcommand\G{\mathcal{G}}

\renewcommand{\O}{\mathcal{O}}

\newcommand{\R}{\mathbb{R}}
\newcommand{\C}{\mathbb{C}}
\renewcommand{\c}{\mathcal{C}}

\newcommand\lie[1]{\mathfrak{#1}}
\renewcommand{\k}{\lie{k}}

\newcommand{\g}{\lie{g}}

\newcommand{\z}{\lie{z}}
\renewcommand{\t}{\lie{t}}
\newcommand{\Alc}{\lie{A}}

\newcommand{\on}{\operatorname}

\newcommand{\length}{\on{length}}

\newcommand{\Ad}{ \on{Ad} } 
\newcommand{\Hol}{ \on{Hol} }

\newcommand{\Vol}{  \on{Vol}}
\newcommand{\diag}{  \on{diag}}

\newcommand\dirac{/\kern-1.2ex\partial} 
\newcommand\qu{/\kern-.7ex/} 
\newcommand{\Waff}{W_{\text{aff}}} 

\usepackage[cmtip,matrix,arrow]{xy}
\CompileMatrices
\SilentMatrices
\SelectTips{cm}{}
\newdir{((}{{\lhook\kern-.1em}}
\newdir{))}{{\rhook\kern-.1em}}
\newdir{ >}{{}*!/-5pt/@{>}}

\newcommand{\hra}{\hookrightarrow}

\newcommand{\bib}{\bibitem}

\renewcommand{\d}{{\mbox{d}}}
\newcommand{\ol}{\overline}

\newcommand\Phinv{\Phi^{-1}}

\newcommand\Sig{\Sigma}
\newcommand\sig{\sigma}

\newcommand\Om{\Omega}
\newcommand\om{\omega}

\newcommand{\f}{\frac}

\newcommand{\p}{\partial}
\renewcommand{\l}{\langle}
\renewcommand{\r}{\rangle}

\newcommand{\ti}{\tilde}

\begin{document}

\sloppy

\title[Moduli spaces of flat connections]
{Moduli spaces of flat connections on 2-manifolds, 
cobordism, and Witten's volume formulas}
\author{E. Meinrenken}
\thanks{Supported by a Feodor Lynen fellowship from the Humboldt foundation.}
\address{Massachusetts Institute of Technology, Department
of Mathematics, Cambridge, Massachusetts 02139}
\email{mein@math.mit.edu}
\author{C. Woodward}
\thanks{Supported by an NSF Postdoctoral Fellowship.}
\address{Harvard University, Department of Mathematics, 1 Oxford Street,
Cambridge, Massachusetts 02138}
\date{May 20, 1997}
\email{ woodward@math.harvard.edu}
\maketitle

\section{Introduction}

According to Atiyah-Bott \cite{AB},\cite{A} the moduli space of flat
connections on a compact oriented 2-manifold with prescribed holonomies around
the boundary is a finite-dimensional symplectic manifold, possibly
singular.  A standard approach \cite{W1,W2,V} to computing invariants
(symplectic volumes, Riemann-Roch numbers, etc.)  of the moduli space
is to study the ``factorization'' of invariants under gluing of 
2-manifolds along boundary components. Given such a factorization 
result, any choice of a ``pants decomposition'' of the 2-manifold 
reduces the computation of invariants to the 
three-holed sphere.

Consider a  compact, connected, simple 
Lie group $G$ with maximal torus $T$, and let $\Alc\subset\t$ be some 
choice of a fundamental Weyl alcove. We interpret $\Alc$ as the set 
of conjugacy classes in $G$ since for every conjugacy class 
$\c\subset G$ there is a unique $\mu\in\Alc$ with $\exp(\mu)\in \c$.
For $\mu_1,\mu_2,\mu_3\in\Alc$ let 
$$ \M(\Sig_0^3,\mu_1,\mu_2,\mu_3)$$
denote the moduli space of flat $G$-connections on the three-holed
sphere $\Sig_0^3$ for which the holonomies around the three 
boundary components lie in the conjugacy classes labeled by the $\mu_j$.

Thanks to the following result of L. Jeffrey \cite{J}, the moduli
spaces $\M(\Sig_0^3,\mu_1,\mu_2,\mu_3)$ are well-understood when the
holonomies are small.  Use the normalized invariant inner product on
$\g$ to identify $\g\cong\g^*$ and let $\O_{\mu_j}$ be the coadjoint
orbit through $\mu_j$, equipped with the Kostant-Kirillov-Souriau
symplectic form, and $\O_{\mu_j}^-$ the same space with the opposite
symplectic form.
\begin{theorem}[Jeffrey] \label{lisajeffrey}
For $\mu_1,\mu_2,\mu_3\in\Alc$ 
sufficiently small, there is a symplectomorphism 
$$ \M(\Sig_0^3,\mu_1,\mu_2,\mu_3)\cong 
\O_{\mu_1}^-\times \O_{\mu_2}^- \times \O_{\mu_3}^-\qu G $$
of the moduli space with a symplectic reduction of a triple product 
of coadjoint orbits.
\end{theorem}
In \cite{MW2} we prove the following result for more general
holonomies:
\begin{theorem} \label{MainTheorem}
For $\mu_1,\mu_2,\mu_3\in \on{int}(\Alc)$ generic  there is an
oriented orbifold cobordism
\begin{equation}
\M(\Sig_0^3,\mu_1,\mu_2,\mu_3)\sim \coprod_{w\in\Waff^+}
(-1)^{\length(w)}\, (\O_{\mu_1}^-\times \O_{\mu_2}^- \times \O_{w\,\mu_3}^-)
\qu G.\label{MainCob}
\end{equation}
Here the signs indicate a change in orientation relative to the
symplectic orientation, and $\Waff^+$ is the set of all 
$w$ in the affine Weyl group $\Waff$ such that 
$w\Alc\subset\t_+$.  The symplectic forms extend 
to a closed 2-form over the cobordism. For $G=\on{SU}(n)$ both 
sides are smooth manifolds and the cobordism is a manifold 
cobordism.
\end{theorem}  
The genericity assumption guarantees that both sides have at worst 
orbifold singularities. Cobordisms of manifolds together with closed 
2-forms were 
introduced and studied by Ginzburg-Guillemin-Karshon \cite{GGK}.  

The main point of this article is not
to give a complete proof of the theorem (which will appear
in \cite{MW2}) but rather to explain the main ideas and to present as 
an application an elementary proof of Witten's formulas for 
symplectic volumes.  These formulas were proved in Witten's 
original paper \cite{W1}; alternative proofs  
appear for $\on{SU}(2)$ in 
\cite{D1,T,JW1,JW2} and for higher rank in \cite{L1,L2,JK}. 
In a forthcoming paper, we use similar cobordism techniques to 
compute the Verlinde numbers.

This article is organized as follows.  In the following section we 
introduce the Guillemin-Ginzburg-Karshon notion of 
cobordism of Hamiltonian spaces and give some examples. 
In section 3 we recall the construction 
of the symplectic form on moduli spaces for 2-manifolds with 
boundary.  In section 4 we explain our construction of the cobordism 
(\ref{MainCob}), and in section 5 we use (\ref{MainCob})
to calculate volumes of moduli spaces. In the appendix we 
collect some background material on Hamiltonian group 
actions and symplectic reduction. 
\section{Hamiltonian cobordism}\label{Cobordism}
Let $K$ be a compact Lie group. By a {\em Hamiltonian $K$-manifold}, we 
mean a triple $(M,\om,\Phi)$ consisting of a manifold $M$ with an 
action of $K$, 
an invariant closed 2-form $\om\in\Om^2(M)^K$, and an 
equivariant {\em moment map}
$\Phi\in C^\infty(M,\k^*)^K$ satisfying 
\begin{equation}
\label{MomentCondition}\d\l\Phi,\xi\r=\om(\xi_M,\cdot)
\end{equation}
for all $\xi\in\k$, where $\xi_M$ denotes the fundamental vector field.

Introducing the complex $\Om_K(M)=(\Om(M)\otimes S(\k^*))^K$ of 
equivariant differential forms, with equivariant differential 
$\ti{d}\alpha(\xi)=\d\alpha(\xi) - 2\pi i\,\iota(\xi_M)\,\alpha(\xi)$, 
Equation
(\ref{MomentCondition})  may be rephrased as the cocycle condition for 
the equivariant 2-form 
$$\ti{\om}=\om+2\pi i\Phi.$$
We denote its equivariant cohomology class by $[\ti{\om}]\in 
H^*_K(M)$. If $\om$ is non-degenerate it is called symplectic. 
In the symplectic case, the 
{\em Liouville volume} is defined as $\Vol(M)=\int_M 
\exp(\om)$. 

For any $\mu\in\k^*$, the reduced space $M_\mu$ is defined as the 
quotient 
$$M_\mu:=\Phinv(\mu)/K_\mu.$$
For $\mu=0$ we will also use the notation $M_0=M\qu K$. If $\mu$ is a
regular value of $\Phi$, the action of the stabilizer group $K_\mu$ on
$\Phinv(\mu)$ is locally free and the reduced space $M_\mu$ is an
orbifold, with a naturally induced closed 2-form $\om_\mu$. Given a
second Hamiltonian action of a Lie group $H$ on $M$ such that the
actions of $H$ and $K$ commute, the $K$-reduced space $M_\mu$ becomes
a Hamiltonian $H$-orbifold. If $\om$ is symplectic then $\om_\mu$ is
symplectic.

The notion of {\em cobordism} of Hamiltonian spaces 
was recently introduced 
by Guillemin-Ginzburg-Karshon \cite{GGK}. 
\begin{definition}\label{DefinitionCobordism} \cite{GGK}
Let $(M_1,\om_1,\Phi_1)$ and $(M_2,\om_2,\Phi_2)$ be 
two oriented Hamiltonian $K$-manifolds with proper moment maps. 
We call $M_1$ and $M_2$ cobordant and write 
$M_1\sim M_2$ if there exists an oriented Hamiltonian 
$K$-manifold with boundary 
$(N,\om,\Phi)$ with proper moment map
such that $\p N=M_1\cup (-M_2)$ and such that 
$\om$ resp. $\Phi$ pull back to $\om_i$ resp. $\Phi_i$. 
\end{definition}

(The properness assumption on the moment map is crucial since
otherwise this notion of cobordism would be more or less void.)  There
is an analogous definition of cobordisms of Hamiltonian orbifolds.

One is mainly interested in the case that the ``ends'' of the cobordism are 
symplectic. Their orientation, however, need not agree 
with the symplectic orientation. Letting $M_{j,r}$ denote 
the connected components of $M_j$, we have as an obvious cobordism 
invariant the {\em signed} symplectic volume: 
$$ \sum_r \pm \Vol( M_{1,r})=
\sum_r \pm \Vol( M_{2,r})
$$
where the signs are $+1$ or $-1$  
according whether or not the orientation agrees 
with the symplectic orientation.

The simplest example of a Hamiltonian cobordism is provided 
by two cohomologous equivariant 2-forms $\ti{\om}_j=\om_j+2\pi i\Phi_j$
on a compact oriented $G$-manifold $M$. Thus suppose     
$M_1=M_2=M$ 
and that there exists a 1-form $\beta\in\Om^1(M)^K$ such that 
$$\om_2-\om_1=\d\beta,\,\, \Phi_2-\Phi_1=-\beta^\sharp$$
where $\beta^\sharp:\,M\to \g^*$ is defined by $\l\beta^\sharp,\xi\r
=\iota(\xi_M)\beta$. 
Let 
$N=M\times [0,1]$ with points denoted $(m,t)$, and let
$$\om=\om_1+\d(t\beta),\,\,\Phi(m,t)=(1-t)\Phi_1(m)+t\Phi_2(m).$$
Then $(N,\om,\Phi)$ provides a cobordism of Hamiltonian $K$-manifolds,
$M_1\sim M_2$.

One of the main results in \cite{GGK} is that any compact symplectic 
Hamiltonian $T$-manifold with $(M,\om,\Phi)$ with isolated fixed points 
is cobordant in the above sense to its ``polarized'' 
linearization around the fixed 
points. More precisely, let $a_j(p)\in\t^*$, $j=1,\ldots,\,
n=\f{1}{2}\dim M$ be the weights for the $T$-action at the fixed points. 
Thus $T_pM\cong\C^n$, and the induced $T$-action has moment map 
$$ \Psi_p:\,T_pM\to \t^*,\,z\mapsto \Phi(p)-\f{1}{2}\sum_j\,|z_j|^2\,a_j(p). $$
Choose a vector $\xi\in\t$ that is not perpendicular to any of the
weights $a_j(p)$. Let
$$ {a}_j^\sharp(p)=\on{sign}\l a_j(p),\xi \r \,a_j(p) $$ 
denote the polarized weights and $\sig_j(p)=\#\{j|\,\l a_j(p),\xi \r<0\}$
the number of sign changes. Let $T_pM^\sharp$ denote the tangent space, with 
the modified $T$-action with moment map 
$$ \Psi_p^\sharp:\,z\mapsto \Phi(p)-\f{1}{2}\sum_j\,|z_j|^2\,a_j^\sharp(p).$$
\begin{theorem}\cite{GGK}\label{GGKtheorem}
 For any compact symplectic Hamiltonian 
$T$-manifold $M$ with isolated fixed points, there is a cobordism of 
Hamiltonian $T$-manifolds
$$ M\sim \coprod_{p\in M^T} (-1)^{\sig(p)}\,T_pM^\sharp.$$
\end{theorem}
If $M_1\sim M_2$ is a cobordism of Hamiltonian $K$-manifolds, and
$\nu\in\k^*$ a regular value of the moment map for the cobordism,
reduction gives rise to a cobordism $(M_1)_\nu\,\sim (M_2)_\nu$.
Using a simple perturbation argument, \cite{GGK} show that if $K_\nu$
is a maximal torus it is in fact sufficient to assume that $\nu$ is a
regular value for the moment maps $\Phi_1,\Phi_2$.

For example, 
in Theorem \ref{GGKtheorem} we have a cobordism 
$$ M_\nu\sim \coprod_{p\in M^T} (-1)^{\sig(p)}\,(T_pM^\sharp)_\nu $$
As an immediate consequence one obtains the formula of 
Guillemin-Lerman-Sternberg \cite{GLS}
$$ \Vol(M_\nu)=\sum (-1)^{\sig(p)} \f{\#\Gamma_p}{\Vol(T)}\,
   \kappa_p(\Phi(p)-\nu),$$
where $\kappa_p$ is the push-forward of the characteristic measure 
on the positive orthant $\R^n_+$ under the map 
$$ \R^n\to \t^*,\,\,x\mapsto \sum_{j=1}^n x_j\, a_j^\sharp(p),$$
$\Gamma_p$ the generic stabilizer 
for the action on $T_pM^\sharp$, and $\Vol(T)$ the volume with respect 
to an inner product on $\t$ which we also use to identify measures and 
functions.

\begin{example}
Let $T$ be the maximal torus of $K$ and $\t^*_+\subset \t^*\subset \g^*$ 
some choice of a positive Weyl chamber. 
Let $\mu\in\on{int}(\t^*_+)$ and $\O_\mu=K\cdot\mu\cong K/T$ 
the coadjoint orbit through 
$\mu$, equipped with the KKS symplectic form $\sig_\mu\in\Om^2(\O_\mu)$.
The moment map for the $K$-action on $\O_\mu$ is the embedding 
$\iota_\mu:\O_\mu\hra\k^*$. 
Let $\O_\mu^{(T)}$ denote the orbit considered as a $T$-manifold, with 
moment map $\text{pr}_{\t^*}\circ\iota_\mu$. 
The $T$-fixed 
points are just the Weyl conjugates $w\cdot\mu$, and taking the polarizing 
vector $\xi$ in $\on{int}(\t_+)$ the polarized weights are 
the positive roots. Thus  
\begin{equation}\label{ClassicalKostant}
\Vol(\O_\mu^{(T)})_\nu = \f{\# Z(G)}{\Vol(T)}\sum_{w\in W}
(-1)^{\length(w)}\kappa(w\mu-\nu) 
\end{equation}
where $\kappa$ is the measure determined by the positive roots. This formula 
is the classical analogue of the Kostant multiplicity formula.
\end{example}

\begin{example}
Suppose $(M,\om,\Phi)$ 
is a compact Hamiltonian $K$-manifold, and  
$\mu\in\t^*_+$. 
Let $X_1=\O_\mu\times M$ be the product, with diagonal $K$-action 
and $\om_{X_1}=\sig_\mu+\om_M$. Let $X_2$ be the associated bundle 
$$  X_2=K\times_T M. $$
There exists a unique closed 2-form $\om_{X_2}$ on $X_2$ for which 
the $K$-action is Hamiltonian, with moment map $\Phi_{X_2}$, 
and such that the pull-back of
$\om_{X_2}$ resp. $\Phi_{X_2}$ to $M$ is given by $\om$ and 
$\text{pr}_{\t^*}\circ\Phi+\mu$. There is a $K$-equivariant diffeomorphism 
$$\phi:\,X_2=\,K\times_T M\to X_1=\O_\mu \times M ,\,\,
[k,m]\mapsto (k\cdot\mu,k\cdot m).$$
We leave it as an exercise to the reader to verify that the pull-back
$\phi^*\ti{\om}_{X_1}$ of the equivariant symplectic form on $X_1$ 
is equivariantly cohomologous to the equivariant form
$\ti{\om}_{X_2}$. It follows that  for any regular value
$\tau\in\on{int}(\t^*_+)$ of the moment maps 
of both $X_1$ and $X_2$, there is
an oriented orbifold cobordism of reduced 
spaces $(\O_\mu \times M)_\tau\sim (K\times_T M)_\tau$. 
The reduction on the right hand side is easily determined: 
Since $\Phi_{X_2}([k,m])=k\cdot(\mu+\text{pr}_{\t^*}\circ\Phi)$,
the $\tau$-level set condition is $\text{pr}_{\t^*}\circ\Phi=
k^{-1}\tau-\mu$. Since this is contained in $\t^*$ only if 
$k^{-1}$ represents an element  $w$ of the Weyl group, it 
follows that the reduction is a disjoint union of 
$T$-reductions of $M$ at values $w\tau-\mu$.
Denoting by $M^{(T)}$ the space $M$ considered as a Hamiltonian 
$T$-manifold, we find a cobordism
$$ (\O_\mu \times M)_\tau\sim \coprod_{w\in W}(-1)^{\length(w)}
M^{(T)}_{w\tau-\mu}.
$$
The resulting formula 
$$ \Vol(\O_\mu \times M)_\tau=\sum_{w\in W}(-1)^{\length(w)}
\Vol(M^{(T)}_{w\tau-\mu})
$$
may be considered the classical analogue to the Frobenius reciprocity
formula. Specializing to the case where $M=\O_\nu$ is a coadjoint orbit 
($\nu\in\on{int}(\t^*_+)$) and combining with the previous example we find 
a cobordism 
$$ (\O_\mu \times \O_\nu)_\tau\sim \coprod_{w_1\in W}\sum_{w_2\in W}
(-1)^{\length(w_1 w_2)} (\k/\t)_{w_1\mu+w_2\nu-\tau}$$
where $\k/\t$ has symplectic structure given by any $T$-invariant
Hermitian metric.\footnote{Any two such structures are equivariantly
symplectomorphic.}  As a consequence we have for symplectic volumes
the formula
\begin{equation}\label{ClassicalSteinberg} 
\Vol(\O_{\mu}\times\O_{\nu})_\tau= \f{\# Z(G)}{\Vol(T)}\, \sum_{w_1\in
W}\sum_{w_2\in W} (-1)^{\length(w_1 w_2)}\kappa(w_1\mu+w_2\nu-\tau)
\end{equation}
which is the classical analogue of the Steinberg formula. This 
was derived in Guillemin-Prato \cite{GP} by a different method.
\end{example}

\section{Moduli Spaces of Flat Connections}

In this section we recall the construction of moduli spaces of 
flat connections over 
2-manifolds. 
Let $G$ be a simple, compact, connected Lie group and $T\subset G$ its maximal
torus.
The integral lattice $\{\xi\in \t|\,\exp(\xi)=1\}$ will be denoted 
by $\Lambda$. We let the {\em weight lattice} $\Lambda^*\subset\t^*$
be its dual and let  
$\mathfrak{R}^+\subset\Lambda^*$ be a system 
of positive roots of $G$. 
\footnote{Note that we work with {\em real}
weights as in \cite{BrD}. Our main motivation for this convention is
that the ``pre-quantizable'' coadjoint orbits $G\cdot\mu$,
($\mu\in\t^*\subset\g^*$) are precisely those through real weight lattice
points.  By contrast Witten \cite{W1,W2} works with {\em infinitesimal}
roots and weights, these differ by a factor $2\pi i$ and give rise to
different powers of $2\pi$ in some of his formulas.} 
Denote by $\t_+\subset \t$ the positive Weyl chamber determined by 
$\mathfrak{R}_+$ and by
$$\Alc=\{\xi\in \t_+|\,\l\alpha_{max},\xi\r\le 1\}$$
the closed fundamental alcove, 
where $\alpha_{max}\in\mathfrak{R}^+$ is the highest root. We identify 
$\Alc$ with the set $\on{Conj}(G)$ of conjugacy classes of $G$ via
\begin{equation}\label{Identifications}
\Alc\cong \t/\Waff \cong T/W \cong G/\Ad(G)=\on{Conj}(G),
\end{equation}
where $W$ is the Weyl group and  $\Waff=\Lambda\rtimes W$ the 
affine Weyl group.

Consider a compact, connected, oriented 2-manifold $\Sig=\Sig^b_h$ of
genus $h$ with $b$ boundary components $B_1,\ldots,\,B_b$.  Let
$\A(\Sig)\cong\Om^1(\Sig,\g)$ denote the space of connections on the
trivial bundle $\Sig\times G$.  We take these connections to be of a
fixed Sobolev class $r>\f{1}{2}$ and consider $\A(\Sig)$ as a Banach
manifold. The gauge group $\G(\Sig)=\on{Map}(\Sig,G)$ consisting of
maps of Sobolev class $r+1$ is a Banach Lie group, which acts on
$\A(\Sig)$ by $g\cdot A=\Ad_g(A)-\d g\,g^{-1}$.  For any boundary $B_j$
and any connection $A\in\A(\Sig)$, the holonomy $\Hol_{B_j}(A)$ around
$B_j$ is well-defined up to conjugacy; here we take the orientation of
the boundary $\p\Sig$ to be opposite to the orientation induced from $\Sig$.  Given
conjugacy classes $\c_1,\ldots,\c_b\subset G$ labeled by
$\mu_1,\ldots,\mu_b\in \Alc$, let
\begin{equation}
\label{ModSpace}
\M(\Sig,\mu_1,\ldots,\mu_b):=\{A\in \A(\Sig)|\,\on{curv}(A)=0, \,\,
\Hol_{B_j}(A)\in \c_j\}/\G(\Sig) 
\end{equation}
be the moduli space of flat connections with holonomies in $\c_j$. 
(This space is independent of the choice of Sobolev class $r$.)
There is a canonical 
isomorphism between this space and the representation variety 
\begin{equation}\label{RepresentationVariety} 
\Big\{(a,c)\in G^{2h}\times 
\c_1\times\ldots\times \c_b \Big|
\,\prod_{i=1}^{h}[a_{2i-1},a_{2i}]
=\prod_{j=1}^b c_j\Big\}/G
\end{equation} 
given by describing each gauge equivalence class  
in terms of its parallel transport. This description shows that 
$\M(\Sig,\mu_1,\ldots,\mu_b)$ is 
a compact, stratified  space (usually singular). 
It was one of the basic observations in Atiyah-Bott \cite{AB} that the
spaces $\M(\Sig,\mu_1,\ldots,\mu_b)$ carry natural symplectic
structures, and that (\ref{ModSpace}) can in fact be viewed as a
symplectic quotient.  Indeed, $\A(\Sig)$ carries a natural symplectic
form
$$ \om_A(a_1,a_2)=
\int_\Sigma a_1 \stackrel{\cdot}{\wedge} a_2
\,\,\,\,\,\,\,\,\,(a_i\in T_A\A(\Sig)\cong\Omega^1(\Sig,\g)), $$
using the normalized invariant inner product on $\g$
for which $\alpha_{max}$ has length $\sqrt{2}$.
\footnote{This inner product is related to the Killing form by
$\l\xi,\eta\r=-(8\pi^2\,n_G)^{-1}\on{tr}(\on{ad}_\xi\on{ad}_\eta)$,
where $n_G$ is the dual Coxeter number of $G$.} 
The
$\G(\Sig)$-action preserves the symplectic form and is in fact
Hamiltonian, with moment map given by
$$\Psi(A)=(\on{curv}(A),\,A|\p\Sig)\in \Om^2(\Sig,\g)
\oplus\Om^1(\p\Sig,\g),$$ 
that is
$$ 
\l\Psi(A),\xi\r=
\int_\Sigma \on{curv}(A)\cdot\xi
+
\int_{\partial \Sigma} \iota^*\,(A\cdot\xi),
$$
for $\xi\in \Om^0(\Sig,\g)=\on{Lie}(\G(\Sig))$.  
This exhibits $\M(\Sig,\mu_1,\ldots,\mu_b)$ as a
symplectic reduction since the pull-backs $A|B_j$ are determined up to
gauge equivalence by the conjugacy class of their holonomies
$\Hol_{B_j}(A)$. 

In the presence of at least one boundary component, it is 
convenient to perform the above reduction in stages.
Let $\G(\p\Sig)=\on{Map}(\p\Sig,G)$ (maps of Sobolev class
$r+\f{1}{2}$) be the gauge transformations of the boundary. Since $G$
is simply connected, the restriction map $\G(\Sig)\to \G(\p\Sig)$ is
surjective. Its kernel $\G_\p(\Sig)$ consists of gauge transformations
that are the identity on the boundary. The moment map for the action
of $\G_\p(\Sig)$ on $\A(\Sig)$ is just the curvature $A\to
\on{curv}(A)$.  Let $\M(\Sig)$ be the symplectic quotient
$$ \M(\Sig)=\A(\Sig)\qu \G_\p(\Sig)=
\{A\in\A(\Sig)|\,\on{curv}(A)=0\}/\G_\p(\Sig). $$
Donaldson proved in \cite{D1} that if the boundary of $\Sig$ is
non-empty ($b\ge 1$), the space $\M(\Sig)$ is a smooth (!) symplectic
Banach manifold. The residual action of
$\G(\p\Sig)=\G(\Sig)/\G_\p(\Sig)$ is Hamiltonian, with moment map
given by $[A]\to \iota_{\p\Sig}^* A\in\Om^1(\p\Sig,\g)$.

Choose orientation preserving parametrizations $B_j\cong S^1$ of 
the boundary components, thereby identifying the gauge group 
$\G(B_j)$ with the loop group $LG=\on{Map}(S^1,G)$ and $\Om^1(B_j,\g)$ 
with $L\g^*:=\Om^1(S^1,\g)$. Thus $\M(\Sig)$ is an example 
of a {\em Hamiltonian $LG^b$-manifold with proper moment map}  
$\Phi:\,\M(\Sig)\to (L\g^*)^b$. It admits a description in terms of 
holonomies, just as the spaces $\M(\Sig,\mu_1,\ldots,\mu_b)$:
\begin{equation}
\label{LGHolonomies}
\M(\Sig)=\Big\{ (a,c,\zeta) \in 
G^{2g}\times G^{b-1}\times (L\g^*)^b\Big|
\prod_{i=1}^{2g}[a_{2i-1},a_{2i}]= \prod_{i=1}^b \, \Ad_{c_i} \, \Hol(\zeta_i)
\Big\}
\end{equation}
where $c_1 = 1$. In this description the action of 
$g=(g_1,\ldots, g_b)\in LG^b$ is given by
$$ g\cdot a_i=\Ad_{g_1(0)}a_i,\ \ g\cdot c_j=g_1(0)\,c_j g_j(0)^{-1}
,\ \ g\cdot \zeta_j=\Ad_{g_j}\cdot\mu_j-\d g_j\,{g_j}^{-1}
$$
and the moment map is the projection to the $(L\g^*)^b$-factor.

Note that the moment map is equivariant 
with respect to the {\em affine} action of $LG$ on $L\g^*$ which corresponds 
to interpreting $L\g^*$ as connections on $S^1$ and $LG$ as the gauge 
group.\footnote{Passing to the central extension $\widehat{LG}$ 
makes this into a linear action, where $L\g^*$ is identified with 
the affine subspace $L\g^*\times\{1\}\subset \widehat{L\g}^*$.
Thus $\M(\Sig)$ can also be viewed 
as a Hamiltonian $\widehat{LG}^b$-manifold where the central circle(s) acts 
with constant moment map $1$.} The basic examples are as follows: 

\begin{example} 
\begin{enumerate}
\item The moduli space $\M(\Sig_0^1)$ of the 1-holed sphere is the based 
      loop group $\Om G=LG/G$. The moment map sends $h\in\Om G$ 
      to $h\cdot 0=-h^{-1}\d h\in L\g^*$.    
\item The moduli space  $\M(\Sig_0^2)$ of the 2-holed sphere is 
      $LG\times L\g^*$. The moment map is $(h,\zeta)\to (h\cdot\zeta,-\zeta)$.
\end{enumerate}
\end{example} 
The Hamiltonian loop group spaces $\M(\Sig)$ have been studied from
various points of view by Donaldson \cite{D1}, Segal \cite{S},
S. Chang \cite{C} and others. They satisfy the following factorization
property, which we learned from S. Martin \cite{M}:
\begin{theorem}\label{GluingEqualsReduction}
Let $\Sig$ be obtained from a possibly disconnected 2-manifold
$\hat{\Sig}$ by gluing two boundary components $B_\pm\subset\p
\hat{\Sig}$. Let $LG$ act on $\M(\hat{\Sig})$ by the diagonal action
corresponding to these two boundary circles; here we reverse the
orientation on $B_-$ so that the $B_-$-moment map is $[\hat{A}]\to
-\iota_{B_-}^*\hat{A}$.  The moduli space $\M(\Sig)$ is given by
symplectic reduction
$$\M(\Sig) = \M(\hat{\Sig}) \qu LG.$$
\end{theorem}
This follows formally because the zero 
level condition is that the restrictions  
$\hat{A}|\,B_\pm$ match up. 
A more detailed proof is given in \cite{MW1}.

We view $\g$ as a subset of $L\g^*$ by
identifying $\mu\in \g$ with the connection $\mu\f{d
\phi}{2\pi}\in\Om^1(S^1,\g)$, so that the diagram  
$$
\xymatrix{{\g} \ar[d] \ar[dr]^{\exp}\\ L\g^* \ar[r]_{\on{Hol}} & G}
$$
commutes.
The moduli spaces $\M(\Sig_h^b,\mu_1,\ldots,\mu_b)$ with fixed holonomies 
are obtained from the Hamiltonian $LG^b$-manifold 
$\M(\Sig_h^b)$ as symplectic reductions: Given 
$\mu_j\in\Alc\subset\t\subset\g\subset L\g^*$ we have
$$ \M(\Sig_h^b,\mu_1,\ldots,\mu_b)=
\M(\Sig_h^b)_{\mu_1,\ldots,\mu_b}=
\Phinv(\mu_1,\ldots,\mu_b)/
(LG)_{\mu_1}\times\ldots\times (LG)_{\mu_b}.$$
To achieve a better understanding of these spaces we need to recall 
a few elementary facts about the affine coadjoint action of $LG$ on 
$L\g^*$. The fact that every 
connection on $S^1$ is
determined up to gauge equivalence by the conjugacy class of its
holonomy means that
$$ L\g^*/LG=G/\Ad(G)=\Alc.$$
That is, every coadjoint $LG$-orbit passes through exactly one point
of the fundamental alcove $\Alc$. For any $\mu\in\Alc$, 
the stabilizer group $(LG)_\mu$ is compact, and depends 
only on the open face $\sig$ of $\Alc$ containing $\sig$. 
It is isomorphic (via the map $g\mapsto g(1)$) to the centralizer 
in $G$ of $\exp(\mu)$. 

The affine coadjoint $LG$-action admits finite 
dimensional slices which can be constructed as follows.
Let $\Alc_\sig$ be the union of all open faces  
$\tau$ of $\Alc$ that contain $\sig$ in their closure; thus $(LG)_\nu
\supset (LG)_\mu$ for $\nu\in\Alc_\sig$. Then $U_\sig:=(LG)_\mu\,
\cdot \Alc_\sig$ is a slice for the action at $\mu$. 

\begin{proposition}
Let $\sig_j\subset \Alc$ be open faces containing $\mu_j$.  
The pre-image 
\begin{equation}\label{cross-section}
\Phinv(U_{\sig_1}\times \ldots \times U_{\sig_b})\subset 
\M(\Sig_h^b)
\end{equation}
is a smooth, finite dimensional symplectic submanifold, and 
is a Hamiltonian $(LG)_{\mu_1}\times \ldots \times 
(LG)_{\mu_b}$-manifold 
with the restriction of $\Phi$ 
serving as a moment map. The 
moduli space (\ref{ModSpace}) is obtained from (\ref{cross-section}) 
as a finite dimensional reduction at $\mu_1,\ldots,\mu_b$.
\end{proposition}

This follows from the symplectic cross-section theorem for Hamiltonian 
$LG$-manifolds with proper moment maps. (See e.g. \cite{MW1} for a 
proof.) Note that the restriction of $\Phi$ is only 
affine-equivariant if at least one of the $\mu_j$ is contained 
in the closed face of $\Alc$ opposite to $\{0\}$. Shifting by 
the vector $(\mu_1,\ldots,\mu_b)\in(L\g^*)^b$ makes it into an 
equivariant moment map. 
The finite dimensional spaces (\ref{cross-section}) are generalized 
versions of  the {\em extended moduli spaces} introduced by 
Chang, Huebschmann and Jeffrey. 

The fact that the moduli spaces (\ref{ModSpace}) can be obtained as
finite-dimensional reductions has a number of implications. 
First of all, it shows that they are stratified
symplectic spaces in the sense of Sjamaar-Lerman \cite{SL} (see
Appendix \ref{Stratified}), and therefore have well-defined symplectic
volumes. We also know that the function
$$(\mu_1,\ldots,\mu_b)\mapsto \Vol(\M(\Sig_h^b,\mu_1,\ldots,\mu_b))$$
is continuous and piecewise polynomial over 
$\Phi(\M(\Sig_h^b))\cap \on{int}(\Alc)^b$ (see Appendix \ref{DuiH}).
Another consequence is that  if $(\nu_1,\ldots,\nu_b)\in\Alc^b$ 
is a regular value and $(\mu_1,\ldots,\mu_b)\in \Alc^b$ sufficiently close to 
$(\nu_1,\ldots,\nu_b)$, there is a symplectic fibration 
$$\M(\Sig_h^b,\mu_1,\ldots,\mu_b)\to \M(\Sig_h^b,\nu_1,\ldots,\nu_b)$$ 
with fiber a product of coadjoint orbits, 
$$(LG)_{\nu_1}\cdot(\mu_1-\nu_1)\times \ldots \times
(LG)_{\nu_b}\cdot(\mu_b-\nu_b)
\subset L\g_{\nu_1}^* \times \ldots \times L\g_{\nu_b}^*
$$
(see Appendix \ref{Fibrations}). 
\section{Construction of the Cobordism}
Instead of reducing with respect to all $LG$-factors at once, we can 
also reduce with respect to only some of the factors. 
This is important because 
the action of $(LG)^{b-1}\subset (LG)^b$ on $\M(\Sig_h^b)$ 
is always free, so that the 
corresponding reduction is always regular. 
Reduced spaces near a given value of the moment map 
are fiber bundles over the reduced space at that value. 
This would be automatic for finite-dimensional Hamiltonian $K$-manifolds
(see Appendix \ref{Fibrations}) and is  
shown for Hamiltonian $LG$-manifolds in \cite{MW2}.
 
Let us apply this to the case of the three-holed sphere. We wish to
understand the reductions $\M(\Sig_0^3)_{\mu_1,\mu_2,\cdot}$ for
$\mu_j$ small.  From the holonomy description (\ref{LGHolonomies}), we
have
$$ \M(\Sig_0^3)_{0,0,\cdot} =\Om G =LG/G $$
as a Hamiltonian $LG$-manifold (alternatively this also follows from the 
``Gluing equals Reduction'' principle, Theorem 
\ref{GluingEqualsReduction}).
The holonomy description tells us furthermore that the zero level set 
for the $LG^2\subset LG^3$-action is an associated bundle 
$LG\times_G(G\times G)$. 

Since the coadjoint orbit $(LG)_0\cdot\mu$ through a small value 
$\mu\in\Alc$ is just the usual coadjoint orbit $\O_\mu$ for $G$, 
we obtain the following description of the nearby reductions 
$\M(\Sig_0^3)_{\mu_1,\mu_2,\cdot}$. Let $LG_0=G$ act on the 
product of coadjoint $G$-orbits $\O_{\mu_1}^-\times 
\O_{\mu_2}^-$ by the diagonal action, with moment map 
\begin{equation} 
\O_{\mu_1}^- \times \O_{\mu_2}^-\ni (\alpha_1,\alpha_2)\mapsto
-(\alpha_1+\alpha_2)\in\g^*.
\end{equation} 
The principal bundle $G\to LG\to \Om G$ has a canonical $LG$-invariant 
connection coming from the splitting of Lie algebras $L\g=\g\oplus \Om\g$.  
The minimal coupling form on the associated bundle 
$LG\times_G(\O_{\mu_1}^-\times \O_{\mu_2}^-)$ 
defined by this connection is the unique closed 
2-form that pulls back to the given 2-form on $\O_{\mu_1}^-\times 
\O_{\mu_2}^-$
and for which the $LG$-action is Hamiltonian. The corresponding moment 
map is given by 
$$LG\times_G(\O_{\mu_1}^-\times \O_{\mu_2}^-)\ni 
 [g,\alpha_1,\alpha_2]\mapsto -\,g\cdot(\alpha_1+\alpha_2). $$
For $\mu_j$ sufficiently small the minimal coupling form 
is symplectic (i.e.  weakly non-degenerate).

\begin{proposition}\label{Smallmu}\cite{MW2}
For $\mu_1,\mu_2\in\Alc$ sufficiently small there is an $LG$-equivariant 
symplectomorphism 
\begin{equation} 
\M(\Sig_0^3)_{\mu_1,\mu_2,\cdot}\cong LG\times_G(\O_{\mu_1}^-\times
\O_{\mu_2}^-).
\end{equation} 
\end{proposition}
As a first application we recover Jeffrey's result, Theorem 
\ref{lisajeffrey}: If $\mu_1,\,\mu_2$ 
are small 
we obtain a solution of $-g\cdot(\alpha_1+\alpha_2)=\mu_3\in\Alc$ 
with $\alpha_j\in \O_{\mu_j}$
only
if $g$ is contained in the subgroup of constant loops, $G\subset LG$.
Hence 
$$ \M(\Sig_0^3)_{\mu_1,\mu_2,\mu_3}=
LG\times_G(\O_{\mu_1}^-\times \O_{\mu_2}^-)_{\mu_3}=
(\O_{\mu_1}^-\times \O_{\mu_2}^-)_{\mu_3}
=\O_{\mu_1}^-\times \O_{\mu_2}^-\times \O_{\mu_3}^-\qu G,
$$
q.e.d. 

For general $\mu_j$, Proposition \ref{Smallmu} does not hold -- 
in fact the induced space  $LG\times_G(\O_{\mu_1}^-\times \O_{\mu_2}^-)$
will not even be symplectic if $\mu_1,\mu_2\in\Alc$ get too big. 
However, just as in the finite dimensional setting 
(see Appendix \ref{Fibrations}) we can make the 
following weaker statement: Let 
$\Alc_0:=\{\mu\in\Alc\,|\l\alpha_{max},\mu\r<1\}$ be the 
set of all points in the alcove such that the isotropy group 
$(LG)_\mu$ is contained in $(LG)_0=G$. 
\begin{theorem}\cite{MW2} \label{Cobordism1}
For all $\mu_1,\mu_2\in \Alc_0$, there exists an $LG$-equivariant
diffeomorphism 
$$ \phi:\,\M_I:=\M(\Sig_0^3)_{\mu_1,\mu_2,\cdot}\to 
\M_{II}:=LG\times_G (\O_{\mu_1}^-\times \O_{\mu_2}^-)$$
and an $LG$-invariant 1-form $\beta$ on $\M_I$ such that the 2-forms
$\om_I,\,\om_{II}$ and the moment maps $\Phi_I,\,\Phi_{II}$ are
related by
$$ 
\phi^*\om_{II}=\om_I+\d\beta,\,\,\phi^*\Phi_{II}=\Phi_I-\l\beta,(\cdot)_M\r.
$$
\end{theorem}
Our method of proof in \cite{MW2} is to combine 
Proposition \ref{Smallmu} with a Duistermaat-Heckman type 
result, using that the equivariant cohomology classes on both sides 
vary linearly with $\mu_1,\mu_2$. The slopes of both changes have to
be equal since for small $\mu_j$ the two spaces are symplectomorphic.

As we explained in Section \ref{Cobordism}, the existence of an equivariant 
diffeomorphism such that the equivariant 2-forms are cohomologous 
implies that there is a cobordism of Hamiltonian $LG$-manifolds 
$$\M_I\sim \M_{II}$$
in the sense of 
Ginzburg-Guillemin-Karshon.
Theorem 
\ref{MainTheorem} follows immediately  
by reduction: $(\M_I)_{\mu_3}\sim (\M_{II})_{\mu_3}$. 
We thus need to compute 
reductions at $\mu_3\in\on{int}(\Alc)$ of the Hamiltonian    
$LG$-manifold $LG\times_G (\O_{\mu_1}^-\times \O_{\mu_2}^-)$. The $\mu_3$-level 
set is given by the equation 
$$ -\,g\cdot (\alpha_1+\alpha_2)=\mu_3. $$
Since $G\cdot \t_+=\g$, we may assume $-(\alpha_1+\alpha_2)\in\t_+$.
Then $g^{-1}$ maps $\mu_3\in \on{int}(\Alc)$ to
$-(\alpha_1+\alpha_2)\in \t_+$, which implies that $g\,\cdot \Alc
\subset \t_+$: In other words, $g^{-1}$ represents an element
$w\in\Waff^+$ with $w\mu_3=-(\alpha_1+\alpha_2)$.  This shows that the
reduction at $\mu_3$ is a disjoint union of all reductions
$(\O_{\mu_1}^-\times \O_{\mu_2}^-)_{w\mu_3}$ as $w$ ranges over over
$\Waff^+$. A careful check of orientations gives Theorem
\ref{MainTheorem}. We see that it is in fact sufficient that
$\mu_1,\mu_2\in\Alc_0$ and $\mu_3\in\on{int}(\Alc)$.

The above argument generalizes to the case of a $b$-holed 
sphere with $b\ge 3$: We have 
$$  \M(\Sig_0^b,\mu_1,\ldots,\mu_b)\sim \coprod_{w\in\Waff^+}
(-1)^{\length(w)}\, (\O_{\mu_1}^-\times \ldots\times \O_{\mu_{b-1}}^- 
\times \O_{w\,\mu_b}^-)
\qu G
$$
if $\mu_1,\ldots,\mu_{b-1}\in\Alc_0$ and $\mu_b\in\on{int}(\Alc)$ are 
generic.

\begin{remark}
For the case of $G=\on{SU}(n)$, the cobordism is actually a 
cobordism of manifolds. The reason why no orbifold singularities
appear is as follows. By the holonomy description, every 
stabilizer $K$ for the action of $LG^b$ on $\M(\Sig_h^b)$ 
is isomorphic (by the map $LG\to G,\,g\mapsto g(1)$) to
an intersection $\bigcap_j Z_{g_j}$ of centralizers of 
elements in $G$. For every $k\in K$ and every $g_j$ in this list,
$g_j$ is contained in $Z_k$, so that $Z_{g_j}$ contains the 
center of $Z_k$. Consequently, $K$ can only be discrete if 
every $Z_k$ is semi-simple. However, precisely for 
$G=\on{SU}(n)$ there are no semi-simple centralizers other 
than $G$ itself, so that $K$ has to be equal to $Z(\on{SU}(n))$. 
Since the cobordism is given by the product 
$\M(\Sig_0^3)_{\mu_1,\mu_2}\times [0,1]$, 
it follows that every discrete stabilizer for the $LG$-action 
is $Z(G)$. Hence any reduced space at a regular value is 
smooth. 
\end{remark}

\section{Witten's volume formulas}
In this section we explain how to obtain Witten's volume 
formulas from the cobordism in Theorem 
\ref{MainTheorem}. Since Liouville volumes of 
symplectic orbifolds are clearly invariants under oriented cobordism,
we have the formula 
\begin{equation}\label{Vol1}
 \Vol(\M(\Sig_0^3,\mu_1,\mu_2,\mu_3))=
\sum_{w\in\Waff^+}
(-1)^{\length(w)}\,\Vol(\O_{\mu_1}\times \O_{\mu_2} \times 
\O_{w\,\mu_3}\qu G).
\end{equation}
The volumes on the right hand side are given by the classical analogue 
(\ref{ClassicalSteinberg})
of Steinberg's formula discussed in Section \ref{Cobordism}.
Combining (\ref{ClassicalSteinberg}) with (\ref{Vol1}) gives an 
explicit formula for the volumes  of 
the moduli space $\M(\Sig^3_0)_{\mu_1,\mu_2,\mu_3}$. 
It is convenient to replace
the sum over  
$\Waff^+ =\Waff/W\cong \Lambda$ by a sum over the integral lattice 
$\Lambda\subset \Waff$,  
using the $W$-invariance of the 
function $\nu\mapsto \sum_{w\in W}(-1)^{\length(w)}
\kappa(w\mu-\nu)$:  
\begin{theorem}\label{SteinbergVolume}
For $\mu=(\mu_1,\mu_2,\mu_3)\in\on{int}(\Phi(\M(\Sig_0^3))\cap \Alc^3)$, 
the Liouville volume of the moduli space 
$\M(\Sig_0^3)_{\mu_1,\mu_2,\mu_3}$ is given by the formula
\begin{eqnarray}\label{SteinbergVol}
\lefteqn{
\Vol(\M(\Sig_0^3)_{\mu_1,\mu_2,\mu_3})}\nonumber \\&=&(-1)^{\f{1}{2}\dim G/T}
\f{\# Z(G)}{\Vol(T)}\,
\sum_{l\in\Lambda}\sum_{w_1,\,w_2\in W}
(-1)^{\length(w_1\,w_2)}\kappa(w_1\mu_1+w_2\mu_2+\mu_3+l)
\end{eqnarray}
(where the $w_i$-summation is to be carried out before the summation over 
the integral lattice $\Lambda$). 
\end{theorem}
Note that this holds without regularity assumptions on the $\mu_j$, 
since it is known from basic properties of DH-measures that the 
volume function is continuous over 
$\on{int}(\Phi(\M(\Sig_0^3))\cap \Alc^3)$, even at 
singular values. 

This formula can be re-cast in the following form due to 
Witten. Let $A:\,T\to \R$ be the $W$-anti-invariant function 
$$ A(\exp \xi)=
\prod_{\alpha\in\mathfrak{R}^+}\,2\,\sin(\pi \l\alpha,\xi\r).$$
We label the irreducible $G$-representations by their dominant weights 
$\lambda\in\Lambda^*_+:=\Lambda^*\cap \t_+$ and let 
$\chi_\lambda:\,G\to \C$ denote the character and $d_\lambda=\chi_\lambda(e)$ 
the dimension.  
\begin{theorem} (Witten formula for the three-holed sphere.)
For $\mu=(\mu_1,\mu_2,\mu_3)\in\on{int}(\Phi(\M(\Sig_0^3))\cap \Alc^3)$, 
the volume of the moduli space of the three-holed sphere 
$\M(\Sig_0^3)_{\mu_1,\mu_2,\mu_3}$
is given by
\begin{equation}\label{Witten1}
\Vol(\M(\Sig_0^3)_{\mu_1,\mu_2,\mu_3})=
\# Z(G)  \f{\Vol(G)}{\Vol(T)^3}\prod_{j=1}^3 A(e^{\mu_j})\,
\sum_{\lambda\in\Lambda^*_+}
\f{1}{d_\lambda}\prod_{j=1}^3 \chi_\lambda(e^{\mu_j}).
\end{equation}
Here $\Vol(T)$ and $\Vol(G)$ are the Riemannian volumes with respect 
to the normalized inner product on $\g$. 
\end{theorem} 
\begin{proof}
Note first that the right hand sides of both (\ref{Witten1}) 
and (\ref{SteinbergVol}) define naturally 
$\Waff$-anti-invariant tempered distributions in 
$\mu_1,\mu_2,\mu_3\in\t$. It is sufficient to show 
that we get the same answer if we apply 
the constant coefficient differential operator 
$${\mathcal D}_3=\prod_{\alpha\in\mathfrak{R}^+}
\l\alpha,\f{\p}{\p{\mu_3}}\r.$$
to both expressions. (We use the fact that $\mathcal{D}_3\,u=0$ has no
$W$-anti-invariant solutions in the space of tempered 
distributions.)  The measure $\kappa$
satisfies ${\mathcal D}_3\kappa(\mu_3) =\delta(\mu_3)$. Using the
Poisson summation formula, we obtain:
\begin{eqnarray*}
\f{1}{\# Z(G)} \lefteqn{ {\mathcal D}_3   
\Vol(\M(\Sig_0^3)_{\mu_1,\mu_2,\mu_3}) }\\
&=& 
\f{(-1)^{{\f{1}{2}\dim G/T}}}{\Vol(T)}\,
\sum_{l\in\Lambda}\sum_{w_1,\,w_2\in W}
(-1)^{\length(w_1\,w_2)}\delta (w_1\mu_1+w_2\mu_2+\mu_3+l)\\
&=&  
\f{(-1)^{{\f{1}{2}\dim G/T}}}{\Vol(T)^2}\,
\sum_{\lambda\in\Lambda^*}\sum_{w_1,\,w_2\in W}
(-1)^{\length(w_1\,w_2)}e^{2\pi i \l w_1\mu_1+w_2\mu_2+\mu_3,\lambda\r}.
\end{eqnarray*}

In this sum, weights $\lambda\in\Lambda^*$ which lie on a wall of
some Weyl chamber do not contribute because
$$ \sum_{w\in W}(-1)^{\length(w)}e^{2\pi i \l w\,\nu,\lambda\r}=0$$
if $\lambda$ has a non-trivial stabilizer in $W$. Since 
$\Lambda^*\cap\on{int}(\t_+)=\rho+\Lambda^*_+$, we can rewrite  
the sum as 
\begin{eqnarray*}
&&
\f{(-1)^{{\f{1}{2}\dim G/T}}}{\Vol(T)^2}\,
 \sum_{\lambda\in\Lambda^*_+}\sum_{w_1,\,w_2,\,w\in W}
(-1)^{\length(w_1\,w_2)}\,e^{2\pi i \l w_1\mu_1+w_2\mu_2+\mu_3,w(\lambda+\rho) \r}\\
&=& 
\f{(-1)^{{\f{1}{2}\dim G/T}}}{\Vol(T)^2}\,
\sum_{\lambda\in\Lambda^*_+}\sum_{w_1,\,w_2,\,w_3\in W}
(-1)^{\length(w_1\,w_2)}\,e^{2\pi i \l w_1\mu_1+w_2\mu_2+w_3 \mu_3, \lambda+
\rho \r}\\
&=&{\mathcal D}_3\,  
\f{(-1)^{{\f{1}{2}\dim G/T}}}{\Vol(T)^2}\,
  \sum_{\lambda\in\Lambda^*_+}
\f{(2\pi i)^{-{\f{1}{2}\dim G/T}}}{\prod_{\alpha\in\mathfrak{R}^+}\l\alpha,\lambda+\rho\r}
\prod_{j=1}^3\big(\sum_{w_j \in W} (-1)^{\length(w_j)}e^{2\pi i \l w_j\mu_j,\,\lambda+\rho\r}\big).
\end{eqnarray*}
The last expression is identified with the right hand side of 
(\ref{Witten1}), using the Weyl character formula 
\begin{equation} \chi_\lambda(\exp \xi)=\f{\sum_{w\in W}(-1)^{\length(w)}
e^{2\pi i\l w(\lambda+\rho),\xi\r}}
{i^{\f{1}{2}\dim G/T}\,A(e^\xi)}
 \end{equation}
and the Weyl dimension formula
\begin{equation} d_\lambda=\f{\prod_{\alpha\in\mathfrak{R}^+}\l \alpha,\lambda+\rho\r}
{\prod_{\alpha\in\mathfrak{R}^+}\l \alpha,\rho\r} = (2\pi)^{
\f{1}{2}\dim G/T}\Vol(G/T) \prod_{\alpha\in\mathfrak{R}^+}\l
\alpha,\lambda+\rho\r .\end{equation}
Here $\Vol(G/T)$ is the Riemannian volume 
of $G/T$ with respect to the inner product on $\g$
(see e.g. \cite{BGV}). 
\end{proof}

\begin{remark}
The differential equation for 
${\mathcal D}_3\Vol(\M(\Sig_0^3)_{\mu_1,\mu_2,\mu_3})$
appearing in this proof is interpreted in Jeffrey-Weitsman
\cite{JW2}
in terms of symplectic volumes of intersections 
of divisors in the moduli space. 
\end{remark} 
From the volume formula for the three-holed sphere, formulas for
the general case are obtained by gluing. For this we need: 
\begin{proposition}\label{VolumeFactorization}
Let $\Sig=\Sig_h^b$ be obtained from a possibly disconnected 
2-manifold $\hat{\Sig}$ by gluing two 
boundary components $B_\pm\subset\p\Sig$. Suppose 
$\mu=(\mu_1,\ldots,\mu_b)\in \on{int}\,\Alc^b$ is such that 
$\M(\Sig,\mu_1,\ldots,\mu_b)$ contains at least one connection 
with stabilizer $Z(G)$.
Then  
$$ \Vol(\M(\Sig,\mu_1,\ldots,\mu_b))=
\f{1}{k}\int_{\Alc}
\Vol(\M(\hat{\Sig},\mu_1,\ldots,\mu_b,\nu,*\nu))|\d\nu|
$$
Here the measure $|\d\nu|$ on $\Alc\subset \t$ is the normalized measure 
for which $\t/\Lambda^*$ has measure 1, and $k=1$ if $\hat{\Sig}$ 
is connected and equal to $\# Z(G)$ if $\hat{\Sig}$ is disconnected.  
\end{proposition}

This is proved in Jeffrey-Weitsman \cite{JW2}. It also follows from 
the ``Gluing equals Reduction'' principle, 
Theorem \ref{GluingEqualsReduction} since the proof of 
Theorem 
\ref{VolumeFactorization1} goes through for Hamiltonian loop 
group actions with proper moment maps (\cite{MW1}, Proposition 3.12).

The reason for the factor $\f{1}{k}$ is that
the generic stabilizer for the $LG\times LG$-action 
on $\M(\hat{\Sig},\mu_1,\ldots,\mu_b,\cdot,\cdot)$ is $Z(G)$ if 
$\hat{\Sig}$ is connected, $Z(G)\times Z(G)$ otherwise.   
Carrying out the integrations gives: 
\begin{proposition}
Suppose $2h+b\ge 3$. 
For all $\mu_1,\ldots \mu_b\in\on{int}(\Alc)$ such that 
$\M(\Sig_h^b,\mu_1,\ldots,\mu_b)$ contains the gauge 
equivalence class of at least one connection with 
stablizer $Z(G)$, the Liouville 
volume of the moduli space $\M(\Sig_h^b,\mu_1,\ldots,\mu_b)$ 
is given by the formula 
$$ 
\Vol(\M(\Sig_h^b,\mu_1,\ldots,\mu_b))=
\#\,Z(G)
\f{\Vol(G)^{2h-2+b}}{\Vol(T)^b}
\prod_{j=1}^b A(e^{\mu_j})\,
\sum_{\lambda\in\Lambda^*_+}
\f{1}{d_\lambda^{2h-2+b}}\prod_{j=1}^b \chi_\lambda(e^{\mu_j}).
$$
\end{proposition}

\begin{proof}
This follows from Proposition \ref{VolumeFactorization} 
and (\ref{Witten1})
by iterated gluing of 
3-holed spheres. We have to use that $A(e^\xi)\,A(e^{*\xi})
=|A(e^\xi)|^2$ and that by Weyl's integration formula
and the orthogonality relations of irreducible characters,
$$ \int_{\Alc} \chi_{\lambda_1}(e^\xi)\chi_{\lambda_2}(e^\xi)\,|A(e^\xi)|^2\, |\d\xi|
=\left\{\begin{array}{cc} 1 & \text{if } \lambda_1=*\lambda_2\\
0 & \text{otherwise}\end{array}\right.$$
where $|\d\xi|$ is the measure on $\t$ which is normalized
with respect to the lattice $\Lambda$. In Proposition 
\ref{VolumeFactorization} we need to use instead Lebesgue measure 
normalized with respect to $\Lambda^*$. Since these two measures 
differ by a factor $\Vol(T)^2$ (where $\Vol(T)$ is the volume of $T$ with 
respect to the inner product on $\t\subset \g$), 
we get an overall factor 
$\Vol(T)^{2(3h-3+b)}$ from the $3h-3+b$ gluing circles, which combines 
with the factor $ \Vol(T)^{-3(2h-2+b)}$ corresponding to the number 
of 3-holed spheres to a factor $\Vol(T)^{-b}$.
\end{proof}

%
%

\begin{remark}
The moduli space $\M(\Sig_h^0)$ for a surface without boundary and
genus $h\ge 2$ always contains a connection with stabilizer equal to
$Z(G)$.
\end{remark}

Finally, we would like to get rid of the assumption that the $\mu_j$ 
lie in the interior of the fundamental alcove.
For $\nu=(\nu_1,\ldots,\nu_b)\in\Alc^b$ such that $\Phinv(\nu)$
contains a connection with stabilizer $Z(G)$,
we have by Appendix \ref{Stratified}, Equation (\ref{LimitingFormula}) 
$$ \Vol\,\M(\Sig_h^b,\nu_1,\ldots,\nu_b)=\lim_{\stackrel{\mu_j\to\nu_j}
{\mu_j\in\on{int}(\Alc)}}
\f{\Vol\,\M(\Sig_h^b,\mu_1,\ldots,\mu_b)}{\Vol(K_1\cdot (\mu_1-\nu_1))\ldots
\Vol(K_b\cdot (\mu_b-\nu_b) )} $$ 
where $K_j=(LG)_{\nu_j}$. Note that $K_j$ contains the maximal torus 
$T$ of $G$. A natural choice of positive Weyl chamber for $K_j$ is 
given by $\t^*_{+,j}=\R_+\cdot (\Alc-\nu_j)$, and the corresponding 
set $\mathfrak{R}_j^+$ of positive roots of $K_j$ is the set of all 
roots $\alpha$ of $G$ 
such that $\l\alpha,\nu_j\r\in\{0,-1\}$.

The Liouville volume of a coadjoint orbit $K_j\cdot (\mu_j-\nu_j)$ is 
related to the Riemannian volume of $K_j/T$ by 
\begin{equation}\label{VolCoad}
\Vol(K_j\cdot (\mu_j-\nu_j))=
(2\pi)^{\f{1}{2}\dim K_j/T} \prod_{\alpha\in \mathfrak{R}_j^+ }\l \alpha,\mu_j-\nu_j\r\,\,
\Vol(K_j/T)
\end{equation}
(see e.g. \cite{BGV}). Since $K_j\cong Z_{\exp(\mu_j)}$  
the factor $\Vol(K_j/T)$ combines with $\Vol(G/T)^{-1}$
to give the Riemannian volume of the conjugacy class $\mathcal{C}_{\mu_j}
=G/Z_{\exp(\mu_j)}$. Moreover
$$ \lim_{\mu_j \to\nu_j}
\f{ \prod_{\alpha\in \mathfrak{R}^+} 2\sin
(\pi\l\alpha,\mu_j\r)        }
{\prod_{\alpha   \in \mathfrak{R}_j^+  }\,2\pi \l \alpha,\mu_j-\nu_j\r}=
\prod_{\stackrel{\alpha\in \mathfrak{R}^+}{\l\alpha,\nu_j\r\not\in\{0,1\}} }
2\sin(\pi\l\alpha,\nu_j\r)  
$$
so that we  obtain: 
\begin{theorem}[Witten Formula]
Suppose $2h+b\ge 3$. Let $\mu=(\mu_1 ,\ldots,\mu_b)\in\Alc^b$  be 
such that the level set $\Phinv(\mu)$ 
contains a connection with stabilizer $Z(G)$. 
The volume of the moduli space of the 
2-manifold $\Sig^b_h$ with fixed 
holonomies $\mu_1,\ldots,\mu_b$ is given by the formula 
\begin{eqnarray*}
\lefteqn{ 
\Vol(\M(\Sig^b_h,\mu_1,\ldots,\mu_b))} \\&=&
\# Z(G)\, \Vol(G)^{2h-2}
\prod_{j=1}^b \big(\Vol(\mathcal{C}_{\mu_j} )\!\!
\prod_{\stackrel{\alpha\in\mathfrak{R}^{+}}{\l\alpha,\mu_j\r\not\in\{0,1\}   }}
\!\!2\sin (\pi \l\alpha,\mu_j\r) \big)
\sum_{\lambda\in\Lambda^*_+}\f{1}{d_\lambda^{2h-2+b}}
\prod_{j=1}^b  \chi_\lambda(e^{\mu_j})\,
\end{eqnarray*}
In particular, 
$$ 
\Vol(\M(\Sig_h^0))=\#\,Z(G)
\Vol(G)^{2h-2}
\sum_{\lambda\in\Lambda^*_+}
\f{1}{d_\lambda^{2h-2}}
$$
\end{theorem}

The above formulas can be made more explicit for $G=SU(2)$.  
Viewed as an element of $\t$, the exponential of $\rho$
is the diagonal matrix $-I\in Z(\on{SU}(2))$ and 
the fundamental alcove is the interval $[0,\rho]$.
Since $\rho=\f{\alpha}{2}$ has length $\f{1}{\sqrt{2}}$,  
we have $\Vol(T)=\sqrt{2}$. Also, $\Vol(\on{SU}(2)/T)=\f{1}{2\pi}$ 
and $\#Z(G)=2$. 
The character of the $k$-dimensional representation is given 
by 
$$\chi_k(\exp(t\rho))=\f{\sin(\pi k t)}{\sin(\pi t)}.$$
For $b=0$ the above formula reads
$$ \Vol(\M(\Sig_h^0))=\f{2^h}{(2\pi)^{2h-2}}\zeta(2h-2). $$
In particular, $\Vol(\M(\Sig_2^0))=\f{1}{\pi^2}\zeta(2)=\f{1}{6}$
which matches with the well-known fact $\M(\Sig_2^0)\cong \C P(3)$.
For $b=1$, $\mu=t\rho$ we have to distinguish the three cases $0<t<1$ and 
$t=0,1$. For $t=0$ we have $\Vol(\M(\Sig_h^1,0))=\Vol(\M(\Sig_h^0))$ as 
expected. For $t=1$ we find $\chi_k(\exp \rho)={(-1)^{k+1}}\,{k}$ and therefore
$$\Vol(\M(\Sig_h^1,\rho))=\f{2^h}{(2\pi)^{2h-2}}
\sum_{k=1}^\infty \f{(-1)^{k+1}}{k^{2h-2}}
=\f{2^h-2^{2-h}}{(2\pi)^{2h-2}}\zeta(2h-2).$$
For $0<t<1$ we have 
$$ \Vol(\M(\Sig_h^1,t\rho))=
2^h \sum_{k=1}^\infty
\f{1}{(2\pi k)^{2h-1}}\sin(\pi k t).$$ 
The Fourier summation gives Bernoulli polynomials, see \cite{D2}.
In particular, for $h=1$ the result is the sawtooth function
$$ \Vol(\M(\Sig_1^1,t\rho))=2(\f{t}{2}-[\f{t}{2}]). $$

\begin{appendix}
\section{Background on Symplectic Reduction}
In this appendix we summarize some basic facts about 
reduction of Hamiltonian $K$-manifolds, where $K$ is a 
compact connected Lie group. 
Let $T\subset K$ be the maximal torus of $K$, and $\t^*_+\subset\t^*
\subset\k^*$ a closed fundamental Weyl chamber. Thus $\t^*_+=\k^*/K$
parametrizes the set of coadjoint orbits in $\k^*$.
The stabilizer $K_\mu$ of 
a point $\mu\in\t^*_+$ depends only on the open face of $\t^*_+$  
containing $\mu$; in particular $K_\mu=T$ for $\mu\in\on{int}(\t^*_+)$.

\subsection{Singular reduced spaces}\label{Stratified}
We start with a brief discussion of singular symplectic quotients. 
Most of this material is due to Sjamaar-Lerman \cite{SL}.  
Suppose that
$(M,\om,\Phi)$ is a connected symplectic Hamiltonian $K$-manifold
($\dim M=2n$) with proper moment map.  Decomposing $M$ according to
conjugacy classes of stabilizer groups, and then decomposing further
into connected components gives the {\em orbit type stratification}
$$M=\cup M_\lambda.$$
There exists a unique {\em principal stratum}
$M^{prin}$ which is open and dense; the corresponding stabilizer group
$\Gamma\subset K$ (defined up to conjugacy) is called the principal
stabilizer. 
Suppose now that $\mu\in\k^*$ is a possibly singular value 
of $\Phi$. Using local normal forms one can show that
every intersection $\Phinv(\mu)\cap M_\lambda$ is 
smooth, and the decomposition
$$\Phinv(\mu)=\cup_\lambda \Phinv(\mu)\cap M_\lambda$$
makes $\Phinv(\mu)$ into a stratified singular space. 
Moreover, taking the quotient by $K_\mu$ one obtains a stratification 
of the reduced space $M_\mu=\Phinv(\mu)/K_\mu$ 
all of whose strata are symplectic manifolds. It is 
shown in \cite{SL} that every $M_\mu$ has a unique open, dense connected 
principal stratum $M_\mu^{prin}$. If $\Phinv(\mu)$ meets $M^{prin}$ 
then 
$$M_\mu^{prin}=(\Phinv(\mu)\cap M^{prin})/K_\mu.$$
The local structure of the singular space $M_\mu$ is that of an 
iterated symplectic cone (see \cite{SL}, section 6), which implies 
in particular that the singular strata have finite Liouville volume.
One defines $\Vol(M_\mu)$ to be the Liouville volume of $M_\mu^{prin}$. 

Let us now assume that the generic stabilizer $\Gamma$ is discrete,
which implies that $\Phi|\,M^{prin}$ is a submersion. 
From the Liouville form $\om^n/n!$ and any translation invariant 
volume form on $\k^*$, 
the level sets $\Phinv(\mu)\cap M^{prin}$ acquire natural volume forms
and orientations.  
Their volume $\Vol(\Phinv(\mu)):=\Vol(\Phinv(\mu)\cap M^{prin})$ 
is related to the volume of the 
reduced space $\Vol(M_\mu)$ by
\begin{equation}\label{vol}
\Vol(\Phinv(\mu))=\f{\Vol(M_\mu)}{\Vol K\cdot\mu}\,\f{\Vol(K)}
{\#\Gamma}
\end{equation}
where $\Vol(K\cdot\mu)$ is the symplectic volume of the coadjoint
orbit through $\mu$, and $\Vol(K)$ the volume of $K$ with respect to
the dual measure on $\k$.  The function $\mu\mapsto \Vol(\Phinv(\mu))$
is continuous over $\Phi(M^{prin})$.  Combining this with
$\Vol(K\cdot\mu)/\Vol(K\cdot\nu)=\Vol(K_\nu\cdot(\mu-\nu))$,
the symplectic volume of the coadjoint $K_\nu$-orbit through 
$\mu-\nu\in \k_\nu^*$, it follows that for any 
$\nu\in\Phi(M^{prin})\cap\t^*_+$,
\begin{equation}\label{LimitingFormula} 
\Vol(M_\nu)=\lim_{\mu\to\nu}
\f{\Vol(M_\mu)}{\Vol(K_\nu\cdot(\mu-\nu))}.
\end{equation}
This is particularly useful since it allows to compute volumes of 
singular reduced spaces as limits of volumes of regular reduced 
spaces. We emphasize that (\ref{LimitingFormula}) does not hold if 
$\nu\in\t^*_+$ is not contained in $\Phi(M^{prin})$.

\subsection{Duistermaat-Heckman}\label{DuiH}
Suppose that $H$ is another compact Lie group and 
that $(M,\om,(\Phi,\Psi))$ is a compact Hamiltonian $K\times H$-manifold. 
Let 
$\sig$ be an open face of $\t^*_+$. The tangent space to $\sig$ is 
given by $\z(K_\sig)^*$, the dual of the Lie algebra of the 
center of $K_\sig$. Let $\Phi(M)_{reg}$ be the set of 
regular values of $\Phi$.
\begin{theorem}[Duistermaat-Heckman] \label{DH}
Let $M$ be a compact Hamiltonian 
$K\times H$-manifold with moment map $(\Phi,\Psi)$. 
The $K$-reduced spaces $M_\mu$ 
for $\mu$ in a connected component of the set 
$\sig\cap \Phi(M)_{reg}$  
are $H$-equivariantly diffeomorphic, 
and the cohomology class of the $H$-equivariant 2-forms 
$\ti{\om}_\mu$ varies linearly: That is, 
$$ [\ti{\om}_{\mu_1}]-[\ti{\om}_{\mu_2}]=\l\mu_1-\mu_2,\ti{c}\r $$
for a fixed class $\ti{c}\in H^*_H(M_{\mu_1})\otimes \z(K_\sig)$. 
\end{theorem}

An immediate consequence of the DH-theorem is that 
if $(M,\om,\Phi)$ is a compact {\em symplectic} 
Hamiltonian $K$-manifold, the function 
$\mu\mapsto\Vol(M_\mu)$ 
is given by a polynomial on every connected component of 
$\sig\cap \Phi(M)_{reg}$. 

Note that in Theorem \ref{DH} we did not require that $\om$ be 
symplectic. The result is a consequence of two facts in equivariant 
cohomology: First, 
for any compact Hamiltonian $ K\times H$-manifold 
$(M,\om,(\Phi,\Psi))$, with $0$ a regular value of $\Phi$, the 
equivariant cohomology class of $\ti{\om}_0$ is the image of 
the $K\times H$-equivariant class of $\ti{\om}$ under the map 
$$ H_{K\times H}^*(M)\to H^*_{K\times H}(\Phinv(0))\to H^*_H(M\qu K).$$
Second, for any regular value $\mu$ of $\Phi$, zero is a regular value 
for the diagonal action on $M\times \O_\mu^-$, and 
$M_\mu=M\times O_\mu^-\qu K$. Clearly, as $\mu$ varies in $\sig$ 
the equivariant cohomology class of the equivariant 2-form on 
$M\times \O_\mu^-\cong M\times K/K_\sig$ varies linearly.
\subsection{Symplectic Fibrations}\label{Fibrations}
Suppose that $(M,\om,(\Phi,\Psi))$ is a compact 
a Hamiltonian $K\times H$-manifold.
Let $\sig,\tau$ be open faces of $\t^*_+$ with $\sig\subset
\ol{\tau}$, and let  $\nu\in\sig$ and $\mu\in\tau$ regular values 
of $\Phi$, 
contained in the same connected component of $\Phi(M)_{reg}$. 
Choosing an $H$-equivariant diffeomorphism 
$\Phinv(\mu)\cong \Phinv(\nu)$ makes  
$M_\mu$ into a fiber bundle over $M_\nu$, with fiber the coadjoint 
orbit $K_\nu/K_\mu\cong K_\nu\cdot(\mu-\nu)\subset \k_\nu^*$: 
$$ M_\mu\cong \Phinv(\nu)\times_{K_\nu}\,(K_\nu\cdot(\mu-\nu))^-.$$
Let $\theta\in\Om^1(\Phinv(\nu),\k_\nu)^{K_\nu}$ 
be an $H$-invariant principal connection, and let $\sig_{\nu,\mu}$ 
be the KKS-form on $K_\nu\cdot(\mu-\nu)$. The closed 2-form 
on $\Phinv(\nu)\times K_\nu\cdot(\mu-\nu)$ given by 
$$ \iota_\nu^*\om-\sig_{\nu,\mu}-\d\l\alpha,\theta\r$$
(where $\iota_\nu:\Phinv(\nu)\to M$ and $\alpha:\,K_\nu\cdot(\mu-\nu)
\to \k_\nu^*$ are the embeddings) is basic, and therefore descends 
to a closed 2-form on the associated bundle 
$\Phinv(\nu)\times_{K_\nu}\,(K_\nu\cdot(\mu-\nu))^-$
known as the {\em minimal coupling form} of Sternberg \cite{S}. It 
is easy to see that the $H$-action on this bundle is Hamiltonian, with 
moment map naturally induced from $\Psi$.
\begin{theorem}[Fibrations of Reduced Spaces]
\label{FibrationResult}
Let the associated bundle 
\begin{equation}\label{AssociatedBundle}
\Phinv(\nu)\times_{K_\nu}(K_\nu\cdot(\mu-\nu))^-\to M_\nu 
\end{equation}
be equipped with the minimal coupling form for an 
$H$-invariant  principal connection 
$\theta\in\Om^1(\Phinv(\nu),\k_\nu)$. 
The $H$-equivariant closed 2-form on $M_\mu$ is equivariantly cohomologous 
to the equivariant minimal coupling form. If $\om$ is symplectic and 
$\mu$ sufficiently close to $\nu$, the minimal coupling form is 
symplectic, and $M_\mu$ is equivariantly symplectomorphic the 
associated bundle $\Phinv(\nu)\times_{K_\nu}\,(K_\nu\cdot(\mu-\nu))^-$. 
\end{theorem}
Symplectic fibrations by coadjoint orbits are studied in great 
detail in the book \cite{GLS}.

\subsection{Diagonal reduction}
Suppose now that $(M,\om,(\Phi_+,\Phi_-))$ is a compact
connected  symplectic Hamiltonian $K\times K$-manifold. 
(Typically $M$ is the direct product of two Hamiltonian
$K$-manifolds $(M_\pm,\om_\pm,\Phi_\pm)$.)
Suppose also that the generic stabilizers
$\Gamma_{diag}\subset \diag(K)$ for the diagonal $K$-action
on $\Phi_{diag}^{-1}(0)$ and $\Gamma\subset K\times K$ for the 
$K\times K$-action on $M$ are discrete. 
Let 
$*:\,\t^*_+\to \t^*_+$ be the involution defined by
$$ K\cdot(*\mu):=K\cdot (-\mu)\cong (K\cdot\mu)^-.$$  
\begin{theorem}\label{VolumeFactorization1}
The symplectic volume of $M\qu \diag(K)$ is related to the volumes 
of the $K\times K$-reduced spaces by 
$$
\Vol( M\qu\diag(K))=\Vol(T)\f{\# \Gamma}{\# \Gamma_{diag}}\,
\int_{\t^*_+}
\Vol(M_{\mu,*\mu})\,|\d\mu| . 
$$
Here $|\d\mu|$ is any choice of Lebesgue measure on $\t^*$ and 
$\Vol(T)$ the volume of the torus with respect to the induced 
measure on $\t$.  
\end{theorem}
Notice that this result is more or less obvious 
if $K=T$ is a torus. 
The general case reduces to the abelian case as follows: By the 
Guillemin-Sternberg symplectic cross-section theorem, the subset 
$Y=(\Phi_+,\Phi_-)^{-1}(\on{int}(\t^*_+)\times -\on{int}(\t^*_+))$ 
is a smooth symplectic submanifold of $M$, and is a Hamiltonian 
$T\times T$-manifold with the restriction of $(\Phi_+,\Phi_-)$ as a moment 
map. Moreover, $(K\times K)\cdot Y$ is open and dense in $M$ which implies 
that $Y\qu\diag(T)\subset M\qu \diag(K)$ is open and dense. 

\end{appendix}

\end{document}